\documentclass[a4paper,11pt]{article}
\usepackage[utf8]{inputenc}
\usepackage{setspace}
\usepackage{a4}
\usepackage{amssymb,amsmath,amsthm,latexsym}
\usepackage{amsfonts}
\usepackage{amsfonts}
\usepackage{graphicx}
\usepackage{float}
\usepackage{authblk}
\usepackage{amsmath}
\usepackage{xcolor}
\usepackage{lineno,hyperref}
\usepackage{soul}
\usepackage{comment}

\usepackage{mathtools}

 \usepackage{mathtext}
 \restylefloat{table}
\newtheorem{thm}{Theorem}[section]

\newtheorem{pro}[thm]{Proof}

\newtheorem{defn}[thm]{Definition}
\newtheorem{example}[thm]{Example}

\newtheorem{prop}[thm]{Proposition}
\newtheorem{remark}{Remark}
\begin{document}
\begin{center} {\Large \textbf{ AMDS  and  quantum AMDS  Constacyclic codes of length $4p^\varsigma $ over $\mathbb{F}_{{p}^{m}}$ }}
\end{center}

  \begin{center}
    Manasa K J$^{1}$, Divya Acharya$^{2}$, Prasanna Poojary$^{2}$ \\
  $^1$ $[Formaly]$\small{Department of Mathematical and Computational Sciences, National Institute of Technology Karnataka, Surathkal}\\
  $^2$\small{Manipal Institute of Technology, Manipal Academy of Higher Education, Manipal, India}
  {\it \textcolor{blue}{manasakj123@gmail.com;}}{\it \textcolor{blue}{acharyadivya1998@gmail.com; }}\\
		{\it \textcolor{blue}{poojary.prasanna@manipal.edu;}}{\it \textcolor{blue}{poojaryprasanna34@gmail.com; }}\\

 \end{center}

\begin{abstract}

 This paper provides a comprehensive analysis of almost maximum distance separable (AMDS) constacyclic codes of length  $4p^{\varsigma}$  over the finite field  $\mathbb{F}_{p^m}$, where $p$ is an odd prime number.  
 Furthermore, it introduces the construction of quantum AMDS (qAMDS) codes, drawing on the principles of the Calderbank-Shor-Steane (CSS) framework, which enhances their applicability in quantum error correction. This work aims to deepen the understanding of these codes and their potential uses in both classical and quantum computing environments.
\end{abstract}
\textbf{Keywords} Constacyclic codes, Hamming distance, AMDS codes, quantum codes, quantum AMDS codes.
\section{Introduction}

Cyclic and negacyclic codes, which are essential constructs in coding theory, have been the focus of research since the late 1950s, primarily because of their significant practical applications and rich algebraic structures. Among these constructs, constacyclic codes serve as a natural generalization of cyclic codes.
To define a $\beta$ -constacyclic code, we consider a nonzero element $\beta$ from a finite field  $\mathbb{F}_{p^m}$. For a given length $ n$, a $\beta$-constacyclic code is represented as the ideals  $\langle r(x) \rangle$  in the quotient ring  $\frac{\mathbb{F}_{p^m}[x]}{\langle x^n - \beta \rangle}$, where $\langle r(x)\rangle $ is a divisor of $ x^n-\beta.$ 
In the context of these codes, they are categorized into two main types based on their relationship with the characteristic of the finite field. When the code length $n$ is relatively prime to the characteristic of the field, the resulting structure is known as a simple root code.  Otherwise, it is referred  as repeated root codes. The concept of repeated root codes was first explored by Berman $\cite{berman1967semisimple}$ in 1967; this foundational work paved the way for subsequent research by many scholars, including Massey et al. $\cite{massey1973polynomial}$, Roth, and Seroussi $\cite{roth1986cyclic}$. 
However, the most extensive investigations into repeated root codes were carried out by Van Lint $\cite{van1991repeated}$ and later by Castagnoli et al. $\cite{castagnoli1991repeated}$, who significantly contributed to our understanding of their properties and applications. This body of work highlights the importance of constacyclic codes in both theoretical research and practical implementations in coding theory.

Let $ C$ be a code with parameteres $ [n.k.d]_q,$ then $ n, k, d$
must satisfy the singleton bound $ k\leq n-d+1$ \cite{macwilliams1977theory}. A maximum distance separable (MDS) code is one in which $ k=n-d+1.$ For fixed $ n$ and $ k$, an MDS code has the greatest error correcting and detecting capabilities. Thus, constructing MDS codes is important in coding theory. However, designing MDS codes is challenging because meeting the Singleton bound is difficult. Hence, attention shifted to codes that are close approximations to MDS, namely, almost maximum distance separable (AMDS) and near MDS (NMDS) codes. These AMDS and NMDS codes have a connection with BCH codes, Generalized Reed Solomon codes, Ternary Golay codes, and Quadratic residual codes. A linear code with parameters  $ [n,k,d]_q,$ is an AMDS code if $ d=n-k.$ NMDS linear code is one where both the code $ C$ and  its dual $ C^{\perp}$ are AMDS codes. The dual of an MDS code is an MDS code, but this is not true for AMDS codes in general; however, the dual of an NMDS code is guaranteed to be NMDS.
The extensive applications of AMDS and NMDS have been investigated in a series of papers \cite{de1996almost, dodunekov1995near, dodunekova2002almost, faldum1997codes, faldum2002characterization, geng2022class}. These codes have applications in cryptography, secret sharing schemes, design theory, and data storage. The fundamental limit on the parameters of codes in coding theory is the Singleton bound, which is achieved by the MDS codes. Several researchers have worked on MDS codes \cite{dinh2021quantum}, \cite{rani2021quantum}  due to their importance in theory and applications. However, AMDS codes closely approximate the Singleton bound, which has attracted the attention of researchers.

Recently, the development of quantum error-correcting codes with optimal parameters has become a prominent research topic.
To mitigate decoherence and other forms of quantum noise in quantum information, quantum error correction (QEC) was developed.

QEC plays a vital role in quantum computing and quantum communication.
Classical error correcting codes are inadequate for preserving quantum information, hence QEC codes are imperative. The QEC  code was first introduced by Shor \cite{shor1995scheme}  in 1995. Following the groundbreaking research by Calderbank et al. \cite{calderbank1998nested}, the field of quantum information theory has seen rapid growth. Let $ C=[[n,k,d]]_{q}$ be a QEC code. If $ k=n-2(d-1), $ then $ C$ is called a quantum maximum distance separable(qMDS).
 Quantum systems store information in quantum bits (qubits), which can experience various types of errors, including bit flips, phase flips, and depolarizing errors. Quantum Maximum Distance Separable (qMDS) codes serve as a method for protecting quantum states against these errors.
 A quantum extension of AMDS codes is quantum AMDS codes, which play a vital role in protecting data within quantum communication systems.
  These advanced coding techniques are crucial for ensuring the security and integrity of data transmitted through quantum communication systems. By effectively mitigating potential vulnerabilities, quantum AMDS codes significantly enhance the reliability and robustness of quantum data transmission.
   Let $ C=[[n,k,d]]_q$ be  a QEC code with $ k=n-2d,$ then $ C$ is called a qAMDS code. 
   The confidential information remains confidential when quantum states are encoded by an AMDS code.
   qAMDS codes are an effective tool that integrates error correction and security features in quantum information processing.
   In \cite{dinh2024amds}, all AMDS  constacyclic codes and qAMDS codes of length $ p^{\varsigma}$ over the finite field $ \mathbb{F}_{p^m}$ from repeated root codes of prime power length over finite fields are studied.
    In \cite{dinh2025quantum}, \cite{dinh2025amds}
   qAMDS  negacyclic codes of length $ 2p^{\varsigma}$ over $  \mathbb{F}_{p^m}$ and AMDS and qAMDS codes of length $ 3p^{\varsigma}$ over $\mathbb{F}_{p^m} $ studied respectively.    
   The importance of our research is in examining the properties and applications of constacyclic codes of length $ 4p^{\varsigma}$ over the finite field $ \mathbb{F}_{p^m}.$ 
   In our application, we use the CSS (Calderbank-Shor-Steane) construction to systematically identify and derive all quantum AMDS codes. This approach allows us to explore the intricate relationships in quantum coding theory and to uncover a comprehensive set of codes that optimize error-correction capabilities.

The following sections of this paper are structured to provide clear and comprehensive insights: Section 2 contains preliminary results and basic definitions about $ \beta-$ constacyclic codes of length $ 4p^{\varsigma}$ over $  \mathbb{F}_{p^m}.$ In Section 3, we give all AMDS  $ \beta-$ constacyclic codes of length $ 4p^{\varsigma}$ over $  \mathbb{F}_{p^m}.$ In Section 4, we determine quantum AMDS codes constructed from AMDS constacyclic codes of length  $ 4p^{\varsigma}$ over $  \mathbb{F}_{p^m}.$

In Section 5, we will summarize the key findings of our research and draw final conclusions based on the evidence presented throughout this paper. This section will offer insights for future studies.

\section{Preliminaries}
Let $\mathbb{F}_{p^m} $ be a finite field with $ p^{m}$ elements, where $ p$ is an odd prime and $m $ is a positive integer. A nonempty subset  $ C$ of $ \mathbb{F}^{n}_{p^m} $ is called a code of length $ n$ over 
$\mathbb{F} _{{p^m}}. $ If $ C$ is a vector space over $\mathbb{F}_{p^m} $ then it is a linear code.

For an unit $\beta $ of $\mathbb{F}_{{p^m}} $, the $ \beta $ constacyclic shift $ \tau_{\beta}$ on
$\mathbb{F}^{n}_{p^m} $ is $\tau_{\beta}(i_{0},i_{1},...,i_{n-1})=(\beta i_{n-1}, i_{0},....,i_{n-2}). $
If $ \tau_{\beta}(C)=C,$ then $ C $ is a $ \beta-$constacyclic code. 
In particular, for  $ \beta=1$  and $ -1$, such codes are called cyclic and negacyclic codes, respectively.
The following fact is well known\\

Given $ n$ tuples $$ i= (i_{0}, i_{1},...,i_{n-1}),~~ i^{\prime}=( i^{\prime}_{0}, i^{\prime}_{1},....,i^{\prime}_{n-1}    )\in \mathbb{F}^{n}_{p^m} $$ the inner product of $ i,i^{\prime}$ is defined as $ i\cdot i^{\prime}=i_{0}i^{\prime}_{0}+i_{1}i^{\prime}_{1}+\cdots+i_{n-1}i^{\prime}_{n-1}.$
If $ i\cdot i^{\prime}=0$ then $ i$ and $ i^{\prime}$ are  called \textit{orthogonal} codes.
For a linear code $ C $ over $\mathbb{F}_{p^m} ,$  its dual code $ C ^{}$ is given by
$$ C^{\perp}=\{a\in \mathbb{F}_{p^m}| a\cdot b=0 ~\forall~ b\in C\}.$$
 In the following Proposition the dual of a $  \beta-$constacyclic code  is given 
\begin{prop}\cite{dinh2010constacyclic}
	The dual of a $ \beta-$constacyclic code is a $ \beta^{-1}-$ constacyclic code.
\end{prop}

A code $ C$ is called \textit{ dual containing} if $ C^{\perp}\subseteq C.$ If $ C^{\perp}=C$ then the linear code $ C$ is called self dual.

Consider a codeword $ i=( i_{0}, i_{1},....,i_{n-1})\in C .$  The \textit{Hamming weight} of $ i$ denoted by $ wt(i)$, is the count of its nonzero entries. For two codewords $i,j\in C$, the Hamming distance $ d(i,j)$  is the number of positions in which they differ. 
 The Hamming distance $ d(C)$ of a linear code $ C$ is defined as $d(C)=min  \{d(i,j)| i,j\in C, i \neq j \}. $

A complete description of the  structure of $ \beta-$ constacyclic codes of length $ 4 p^{\varsigma }$ over $ \mathbb{F}_{p^{m}}$ appears in \cite{dinh2013repeated}.
We next summarize  the Hamming distances for repeated root  $ \beta-$ constacyclic codes $ C $ of length $ 4 p^{\varsigma }$ over $ \mathbb{F}_{p^{m}}$ obtained by Dinh et al in
\cite{dinh2019hamming};
\begin{thm}\cite{dinh2019hamming}\label{2.2}
	Let $ p$ be an odd prime, $ m$ be a positive integer, $\psi $ and $\phi $ be the nonzero elements of  $ \mathbb{F}_{p^m}.$ Let $ 0\leq {\eta}_{0},\eta_{1}\leq {p-2,}$ and $0\leq{\delta}_{1}\leq{\varsigma-1}$ of the form $ C_{\iota} =  \langle(x^4+\psi)^{\iota}   \rangle$ for $ 0\leq \iota\leq p^{\varsigma}$ and $ x^4+\psi$ be irreducible, then the Hamming distances $ d(C_{\iota})$ are determined by:
	\begin{equation*}
	d(C_{\iota}) =    \begin{cases}
	1, & \text{if } \iota=0,\\
	(\eta_{0}+2)p^{\delta_{0}}, & \text{if } p^{\varsigma }-p^{\varsigma -\delta_{0}}+\eta_{0}p^{\varsigma -\delta_{0}-1}+1
	\leq \iota \leq p^{\varsigma }-p^{\varsigma -\delta_{0}}+(\eta_{0}+1)p^{\varsigma -\delta_{0}-1}
	,\\
	0,& \text{if}~~~ \iota= p^{\varsigma }
	\end{cases}
	\end{equation*}
\end{thm}

\begin{thm}\label{thm:one}\cite{dinh2019hamming}
	If the $ \beta-$ constacyclic codes of length $ 4p^{\varsigma}$ over $ \mathbb{F}_{p^{m}}$ are of the form\\
	$ C_{\iota,\jmath}=\langle(x^2-\phi)^{\iota}(x^2+\phi)^{\jmath}\rangle$ for $ 0\leq \jmath \leq \iota \leq p^{\varsigma},$ then the Hamming distances $ d(C_{\iota \jmath})$ are determined by:
	{\footnotesize \begin{equation*}
	d(C_{\iota \jmath}) =   
    \begin{cases}
	1, & \text{if } \iota = \jmath=0,\\
	2, & \text{if}~ \jmath =0 ~~\text{and}~ 0<\iota \leq p^{\varsigma }\\
	\text{min} \lbrace(\eta_{0}+2)p^{\delta_{0}},2(\eta_{1}+2)p^{\delta_{1}}\rbrace
	& \text{if } p^{\varsigma }-p^{\varsigma -\delta_{0}}+\eta_{0}p^{\varsigma -\delta_{0}-1}+1
	\leq \iota \leq p^{\varsigma }-p^{\varsigma -\delta_{0}}+(\eta_{0}+1)p^{\varsigma -\delta_{0}-1}\\
	& p^{\varsigma }-p^{\varsigma -\delta_{1}}+\eta_{1}p^{\varsigma -\delta_{1}-1}+1
	\leq \jmath\leq p^{\varsigma }-p^{\varsigma -\delta_{1}}+(\eta_{1}+1)p^{\varsigma -\delta_{1}-1}
	,\\
	2(\eta_{1}+2)p^{\delta_{1}}& ~\text{if}~ \iota=p^{\varsigma },\\
	& p^{\varsigma }-p^{\varsigma -\delta_{1}}+\eta_{1}p^{\varsigma -\delta_{1}-1}+1\leq \jmath\leq p^{\varsigma }-p^{\varsigma -\delta_{1}}+(\eta_{1}+1)p^{\varsigma -\delta_{1}-1}\\
	0,& \text{if}~~~ \iota = \jmath= p^{\varsigma }.
	\end{cases}
	\end{equation*}}
	
\end{thm}
\begin{remark}
	By symmetry Hamming distances $ C_{\iota, \jmath }$ for $ \iota \leq \jmath$  is same as in Theorem  \ref{thm:one}. 
	
\end{remark}

\begin{thm} \label{thm:two}\cite{dinh2019hamming}
	Let $ 0\leq \eta_{0},\eta_{1},\eta_{2}\leq p-2 $ and $0\leq \delta_{2}\leq \delta_{1}\leq \delta_{0}\leq \varsigma -1 .  $ Then the Hamming distances of $ C_{\iota, \jmath ,\mu}=\langle(x-1)^{\iota}(x+1)^{\jmath}(x^2+1)^{\mu}\rangle$ for 
	$ 0\leq  \mu\leq \jmath\leq \iota \leq p^{\varsigma} $ are given as follows:
 {\footnotesize \begin{equation*}
		d(C_{\iota, \jmath ,k}) =    \begin{cases}
		1, & \text{if }  \iota = \jmath=\mu=0,\\
		2, & \text{if}~ \mu=0 ~~\text{and}~ 0\leq \jmath \leq \iota\leq p^{\varsigma } \text{( but not  $\iota = \jmath=0$)}\\
		\text{min} \lbrace(\eta_{0}+2)p^{\delta_{0}},2(\eta_{2}+2)p^{\delta_{2}}\rbrace
		& \text{if } p^{\varsigma }-p^{\varsigma -\delta_{0}}+\eta_{0}p^{\varsigma -\delta_{0}-1}+1
		\leq \iota \leq p^{\varsigma }-p^{\varsigma -\delta_{0}}+(\eta_{0}+1)p^{\varsigma -\delta_{0}-1}\\
		& p^{\varsigma }-p^{\varsigma -\delta_{1}}+\eta_{1}p^{\varsigma -\delta_{1}-1}+1
		\leq \jmath\leq p^{\varsigma }-p^{\varsigma -\delta_{1}}+(\eta_{1}+1)p^{\varsigma -\delta_{1}-1},
		\\
		& p^{\varsigma }-p^{\varsigma -\delta_{2}}+\eta_{2}p^{\varsigma -\delta_{2}-1}+1
		\leq \mu \leq p^{\varsigma }-p^{\varsigma -\delta_{2}}+(\eta_{2}+1)p^{\varsigma -\delta_{2}-1},\\
		2(\eta_{2}+2)p^{\delta_{2}}& \text{if} ~~ \iota=p^{\varsigma },\\
		&  p^{\varsigma }-p^{\varsigma -\delta_{1}}+\eta_{1}p^{\varsigma -\delta_{1}-1}+1\leq \jmath\leq p^{\varsigma }\\
		& p^{\varsigma }-p^{\varsigma -\delta_{2}}+\eta_{2}p^{\varsigma -\delta_{2}-1}+1\leq \mu \leq p^{\varsigma }-p^{\varsigma -\delta_{2}}+(\eta_{2}+1)p^{\varsigma -\delta_{2}-1}\\
		0,& \text{if}~~~  \iota = \jmath=\mu= p^{\varsigma }.
		\end{cases}
		\end{equation*}}	
\end{thm}
\begin{remark}
	The Hamming distances of $ C_{\iota, \jmath ,\mu}$ for $\mu \leq \iota \leq \jmath $ is same as in Theorem \ref{thm:two}
\end{remark}
\begin{thm}\label{thm:three}\cite{dinh2019hamming}	
	Let $ 0\leq \eta_{0},\eta_{1},\eta_{2}\leq p-2 $ and $0\leq \delta_{2}\leq \delta_{1}\leq \delta_{0}\leq \varsigma -1 .  $ Then the Hamming distances of $ C_{\iota, \jmath ,\mu}=\langle(x-1)^{\iota}(x+1)^{\jmath}(x^2+1)^{\mu}\rangle$ are given as follows for the case $ 0\leq \jmath\leq \mu \leq \iota \leq p^{\varsigma} $:
 {\footnotesize \begin{equation*}
		d(C_{\iota, \jmath ,\mu}) =    \begin{cases}
		1, & \text{if }    \iota = \jmath=\mu=0,\\
		2, & \text{if}~   \jmath =\mu=0 ~ \text{and}~  0<i<p^{\varsigma},\text{or}~ \jmath =0~\text{and}~ 0\leq \mu \leq \iota \leq p^{\varsigma -1} \\
		3,& \text{if}~  \jmath =0, 0<\mu\leq 2p^{\varsigma -1}~ \text{and}~ p^{\varsigma -1}<\iota \leq2p^{\varsigma -1}\\
		4,&\text{if }   \jmath =0, 0<\mu\leq p^{\varsigma}~\text{and}~ 2p^{\varsigma -1}<\iota \leq p^{\varsigma}\\
		
		\text{min} \lbrace(\eta_{0}+2)p^{\delta_{0}},2(\eta_{1}+2)p^{\delta_{1}},\\4(\eta_{2}+2)p^{\delta_{2}}\rbrace,
		& \text{if } p^{\varsigma }-p^{\varsigma -\delta_{0}}+\eta_{0}p^{\varsigma -\delta_{0}-1}+1
		\leq \iota \leq p^{\varsigma }-p^{\varsigma -\delta_{0}}+(\eta_{0}+1)p^{\varsigma -\delta_{0}-1}\\
		& p^{\varsigma }-p^{\varsigma -\delta_{1}}+\eta_{1}p^{\varsigma -\delta_{1}-1}+1
		\leq\mu \leq p^{\varsigma }-p^{\varsigma -\delta_{1}}+(\eta_{1}+1)p^{\varsigma -\delta_{1}-1},
		\\
		& p^{\varsigma }-p^{\varsigma -\delta_{2}}+\eta_{2}p^{\varsigma -\delta_{2}-1}+1
		\leq \jmath \leq p^{\varsigma }-p^{\varsigma -\delta_{2}}+(\eta_{2}+1)p^{\varsigma -\delta_{2}-1},\\
		\text{min}\lbrace2(\eta_{1}+2)p^{\delta_{1}},4(\eta_{2}+2)p^{\delta_{2}}\rbrace& \text{if} ~~ \iota=p^{\varsigma },\\
		&  p^{\varsigma }-p^{\varsigma -\delta_{1}}+\eta_{1}p^{\varsigma -\delta_{1}-1}+1\leq\mu \leq p^{\varsigma }-p^{\varsigma -\delta_{1}}+(\eta_{1}+1)p^{\varsigma -\delta_{1}-1},\\
		& p^{\varsigma }-p^{\varsigma -\delta_{2}}+\eta_{2}p^{\varsigma -\delta_{2}-1}+1\leq \jmath\leq p^{\varsigma }-p^{\varsigma -\delta_{2}}+(\eta_{2}+1)p^{\varsigma -\delta_{2}-1},\\
		4(\eta_{2}+2)p^{\delta_{2}},~ &\text{if}~~ \iota=\mu=p^{\varsigma },\\
		& p^{\varsigma }-p^{\varsigma -\delta_{2}}+\eta_{2}p^{\varsigma -\delta_{2}-1}+1\leq \jmath\leq p^{\varsigma }-p^{\varsigma -\delta_{2}}+(\eta_{2}+1)p^{\varsigma -\delta_{2}-1},\\ 
		0,& \text{if}~~~ \iota = \jmath =\mu= p^{\varsigma }.
		\end{cases}
		\end{equation*}}
\end{thm}
\begin{remark}
	The Hamming distances of $ C_{\iota, \jmath ,\mu}$  for the cases $ \iota \leq\mu \leq \jmath, \iota \leq \jmath\leq \mu $ and $ \jmath \leq \iota \leq \mu $ are as in Theorem \ref{thm:three}.
\end{remark}
\begin{thm}\label{thm:four}\cite{dinh2019hamming}
	Let $ 0\leq \eta_{0},\eta_{1},\eta_{2},\eta_{3}\leq p-2 $ and $0\leq\delta_{3}\leq \delta_{2}\leq \delta_{1}\leq \delta_{0}\leq \varsigma -1 .  $ Let $ 0 \leq \ell\leq \jmath\leq\mu \leq \iota \leq p^{\varsigma }$, then the codes $ C_{\iota, \jmath ,\mu,\ell }=\langle(x-1)^{\iota}(x+1)^{\jmath}(x-\delta)^{\mu}(x+\delta)^{\ell}\rangle$ have the following Hamming distances:
	
 {\footnotesize \begin{equation*}
    d(C_{\iota, \jmath ,\mu,\ell })=
		\begin{cases}
		1, & \text{if } \iota = \jmath=\mu=\ell =0,\\
		2, & \text{if}~ \jmath =\mu=\ell =0 ~\text{and}~ 0<\iota \leq p^{\varsigma},~\text{or}~~ \ell =0~\text{and}~ 0< \jmath \leq \mu \leq \iota \leq p^{\varsigma -1}, \\
		3,& \text{if}~\jmath =\ell =0, 0<\mu\leq p^{\varsigma }~ \text{and}~ p^{\varsigma -1}<\iota \leq p^{\varsigma },\\
		& \ell =0, 0<\jmath \leq\mu \leq 2p^{\varsigma -1}~ \text{and} ~ p^{\varsigma -1}<\iota \leq 2p^{\varsigma -1},\\
		4,&\text{if}~  \ell =0, 0<\jmath \leq\mu \leq p^{\varsigma} ~\text{and}~ 2p^{\varsigma -1}<\iota \leq p^{\varsigma}\\        
		\text{min} \lbrace(\eta_{0}+2)p^{\delta_{0}},2(\eta_{1}+2)p^{\delta_{1}},\\ 3(\eta_{2}+2)p^{\delta_{2}},4(\eta_{3}+2)p^{\delta_{3}}\rbrace,
		& \text{if } p^{\varsigma }-p^{\varsigma -\delta_{0}}+\eta_{0}p^{\varsigma -\delta_{0}-1}+1
		\leq \iota \leq p^{\varsigma }-p^{\varsigma -\delta_{0}}+(\eta_{0}+1)p^{\varsigma -\delta_{0}-1},\\
		& p^{\varsigma }-p^{\varsigma -\delta_{1}}+\eta_{1}p^{\varsigma -\delta_{1}-1}+1
		\leq \jmath\leq p^{\varsigma }-p^{\varsigma -\delta_{1}}+(\eta_{1}+1)p^{\varsigma -\delta_{1}-1},
		\\
		& p^{\varsigma }-p^{\varsigma -\delta_{2}}+\eta_{2}p^{\varsigma -\delta_{2}-1}+1
		\leq \mu \leq p^{\varsigma }-p^{\varsigma -\delta_{2}}+(\eta_{2}+1)p^{\varsigma -\delta_{2}-1},\\
		& p^{\varsigma }-p^{\varsigma -\delta_{3}}+\eta_{3}p^{\varsigma -\delta_{3}-1}+1
		\leq \ell  \leq p^{\varsigma }-p^{\varsigma -\delta_{3}}+(\eta_{3}+1)p^{\varsigma -\delta_{3}-1},\\
		\text{min}\lbrace2(\eta_{1}+2)p^{\delta_{1}},3(\eta_{2}+2)p^{\delta_{2}}\\ 4(\eta_{3}+2)p^{\delta_{3}}\rbrace& \text{if} ~~ \iota=p^{\varsigma },\\
		&  p^{\varsigma }-p^{\varsigma -\delta_{1}}+\eta_{1}p^{\varsigma -\delta_{1}-1}+1\leq \jmath\leq p^{\varsigma }-p^{\varsigma -\delta_{1}}+(\eta_{1}+1)p^{\varsigma -\delta_{1}-1},\\
		& p^{\varsigma }-p^{\varsigma -\delta_{2}}+\eta_{2}p^{\varsigma -\delta_{2}-1}+1\leq\mu \leq p^{\varsigma }-p^{\varsigma -\delta_{2}}+(\eta_{2}+1)p^{\varsigma -\delta_{2}-1},\\
		& p^{\varsigma }-p^{\varsigma -\delta_{3}}+\eta_{3}p^{\varsigma -\delta_{3}-1}+1
		\leq \ell  \leq p^{\varsigma }-p^{\varsigma -\delta_{3}}+(\eta_{3}+1)p^{\varsigma -\delta_{3}-1},\\ 
		\text{min}\lbrace 3 (\eta_{2}+2)p^{\delta_{2}}, 4(\eta_{2}+2)p^{\delta_{2},}\rbrace &\text{if}~~ \iota = \jmath=p^{\varsigma },\\
		& p^{\varsigma }-p^{\varsigma -\delta_{2}}+\eta_{2}p^{\varsigma -\delta_{2}-1}+1\leq\mu \leq p^{\varsigma }-p^{\varsigma -\delta_{2}}+(\eta_{2}+1)p^{\varsigma -\delta_{2}-1},\\
		& p^{\varsigma }-p^{\varsigma -\delta_{3}}+\eta_{3}p^{\varsigma -\delta_{3}-1}+1
		\leq \ell  \leq p^{\varsigma }-p^{\varsigma -\delta_{3}}+(\eta_{3}+1)p^{\varsigma -\delta_{3}-1},\\
		4(\eta_{3}+2)p^{\delta_{3}},& \text{if}~ \iota = \jmath=\mu=p^{\varsigma},\\
		& p^{\varsigma }-p^{\varsigma -\delta_{3}}+\eta_{3}p^{\varsigma -\delta_{3}-1}+1
		\leq \ell  \leq p^{\varsigma }-p^{\varsigma -\delta_{3}}+(\eta_{3}+1)p^{\varsigma -\delta_{3}-1},\\
		0,& \text{if}~ \iota = \jmath=\mu=\ell = p^{\varsigma }.
		\end{cases}
		\end{equation*}}
\end{thm}
\begin{remark}
	The Hamming distances of $ C_{\iota, \jmath ,\mu,\ell }$ for the cases $ \jmath \leq \ell  \leq \mu \leq \iota , \mu \leq \jmath \leq \ell  \leq \iota , \jmath \leq \mu \leq \ell \leq \iota  , \ell  \leq  \iota \leq \mu \leq \jmath,\mu \leq \iota \leq \ell \leq \jmath , \iota \leq\mu \leq \ell  \leq \jmath, \iota \leq \ell  \leq \mu \leq \jmath, \ell  \leq  \iota \leq \jmath\leq \mu, \iota \leq \jmath \leq \ell  \leq \mu,  \jmath \leq \ell  \leq \iota \leq \mu, \iota \leq \mu \leq \jmath \leq \ell ,\mu \leq \iota \leq \jmath \leq \ell ,  k \leq \jmath \leq \iota \leq \ell  , ~~  and ~~ \jmath \leq \mu \leq \iota \leq \ell  $ is same as in Theorem \ref{thm:four}
\end{remark}

\begin{thm}\label{thm:five}\cite{dinh2019hamming}
	Let $ 0\leq \eta_{0},\eta_{1},\eta_{2},\eta_{3}\leq p-2 $ and $0\leq\delta_{3}\leq \delta_{2}\leq \delta_{1}\leq \delta_{0}\leq \varsigma -1 .  $ Let $ 0 \leq \ell \leq\mu \leq \jmath   \leq \iota \leq p^{\varsigma }$, then the codes $ C_{\iota, \jmath ,\mu,\ell }=\langle(x-1)^{\iota}(x+1)^{\jmath}(x-\delta)^{\mu}(x+\delta)^{\ell}\rangle$ have the following Hamming distances: 
	 {\footnotesize\begin{equation*}
     d(C_{\iota, \jmath ,\mu,\ell })=
		\begin{cases}
		1, & \text{if } \iota = \jmath=\mu=\ell =0,\\
		2, & \text{if}~ \mu=\ell =0~ \text{and}~ 0\leq \jmath \leq \iota\leq p^{\varsigma} \text(but~ not~ \iota = \jmath=0,),\\ &\text{or}~ \ell =0~~\text{and}~~ 0<\mu \leq \jmath\leq \iota \leq p^{\varsigma -1}, \\
		4,&\text{if}~  \ell =0, 0<\mu \leq \jmath \leq p^{\varsigma}~ \text{and}~ p^{\varsigma -1}<\iota \leq p^{\varsigma}\\        
		\text{min} \lbrace(\eta_{0}+2)p^{\delta_{0}},2(\beta_ 
		{2}+2)p^{\delta_{2}},\\ 4(\eta_{3}+2)p^{\delta_{3}}\rbrace,
		& \text{if } p^{\varsigma }-p^{\varsigma -\delta_{0}}+\eta_{0}p^{\varsigma -\delta_{0}-1}+1
		\leq \iota \leq p^{\varsigma }-p^{\varsigma -\delta_{0}}+(\eta_{0}+1)p^{\varsigma -\delta_{0}-1},\\
		& p^{\varsigma }-p^{\varsigma -\delta_{1}}+\eta_{1}p^{\varsigma -\delta_{1}-1}+1
		\leq \jmath\leq p^{\varsigma }-p^{\varsigma -\delta_{1}}+(\eta_{1}+1)p^{\varsigma -\delta_{1}-1},
		\\
		& p^{\varsigma }-p^{\varsigma -\delta_{2}}+\eta_{2}p^{\varsigma -\delta_{2}-1}+1
		\leq \mu \leq p^{\varsigma }-p^{\varsigma -\delta_{2}}+(\eta_{2}+1)p^{\varsigma -\delta_{2}-1},\\
		& p^{\varsigma }-p^{\varsigma -\delta_{3}}+\eta_{3}p^{\varsigma -\delta_{3}-1}+1
		\leq \ell  \leq p^{\varsigma }-p^{\varsigma -\delta_{3}}+(\eta_{3}+1)p^{\varsigma -\delta_{3}-1},\\
		\text{min}\lbrace 2(\eta_{2}+2)p^{\delta_{2}},4(\eta_{3}+2)p^{\delta_{3}}\\ \rbrace& \text{if} ~ \iota=p^{\varsigma },\\
		&  p^{\varsigma }-p^{\varsigma -\delta_{1}}+\eta_{1}p^{\varsigma -\delta_{1}-1}+1\leq \jmath\leq p^{\varsigma }-p^{\varsigma -\delta_{1}}+(\eta_{1}+1)p^{\varsigma -\delta_{1}-1},\\
		& p^{\varsigma }-p^{\varsigma -\delta_{2}}+\eta_{2}p^{\varsigma -\delta_{2}-1}+1\leq\mu \leq p^{\varsigma }-p^{\varsigma -\delta_{2}}+(\eta_{2}+1)p^{\varsigma -\delta_{2}-1},\\
		& p^{\varsigma }-p^{\varsigma -\delta_{3}}+\eta_{3}p^{\varsigma -\delta_{3}-1}+1
		\leq \ell  \leq p^{\varsigma }-p^{\varsigma -\delta_{3}}+(\eta_{3}+1)p^{\varsigma -\delta_{3}-1},\\ 
		4(\eta_{3}+2)p^{\delta_{3},} \text{if}~~  \iota = \jmath=\mu=p^{\varsigma },\\
		& p^{\varsigma }-p^{\varsigma -\delta_{3}}+\eta_{3}p^{\varsigma -\delta_{3}-1}+1\leq \ell \leq p^{\varsigma }-p^{\varsigma -\delta_{3}}+(\eta_{3}+1)p^{\varsigma -\delta_{3}-1},\\
		0,& \text{if}~  \iota = \jmath=\mu=\ell = p^{\varsigma }.
		\end{cases}
		\end{equation*}}
\end{thm}
\begin{remark}
	The Hamming distances of $ C_{\iota, \jmath ,\mu,\ell }$ for the cases 
	$\mu \leq \ell  \leq \jmath \leq \iota, \ell  \leq \mu \leq \iota \leq \jmath,  \mu \leq \ell  \leq \iota \leq \jmath ,  \jmath \leq \iota \leq \ell  \leq \mu, \iota \leq \jmath \leq \ell \leq \mu , \iota \leq \jmath\leq \mu \leq \ell ,~~ and~~~ \jmath \leq \iota \leq \mu \leq \ell ,   $ is same as in Theorem \ref{thm:five}	
\end{remark}
\begin{thm}\cite{dinh2019hamming}
	Assume that $ C_{\iota,\mu }=\langle(x^{2}+\gamma \zeta^{\frac{2p^{m-r}+3p^m-1}{8}}x+ \zeta^{\frac{2p^{m-r}+p^{m}-1}{4}})^{\iota}(x^{2}-\gamma \zeta^{\frac{2p^{m-r}+3p^m-1}{8}}x+ \zeta^{\frac{2p^{m-r}+p^{m}-1}{4}})^{\mu})  \rangle$, where $ 0 \leq\mu \leq \iota \leq p^{\varsigma }$ are integers and $  \gamma \in \mathbb {F}_{p^{m}}$ such that $ \gamma ^{2}=-2$. Let $ 0 \leq  \delta_{0},\delta_{1}\leq \varsigma -1 $ and $   0\leq \eta _{0},\eta_{1}\leq p-2 .$ Then, we have
	 {\footnotesize\begin{equation*}
		d(C_{\iota,\mu }) =    \begin{cases}
		1, & \text{if }~ \iota=\mu=0,\\
		2, & \text{if}~ \mu=0 ~~\text{and}~ 0< \iota \leq p^{\varsigma -1}\\
		3& \text{if}~ \mu=0, p^{\varsigma -1}< \iota \leq p^{\varsigma }\\
		
		\text{min} \lbrace(\eta_{0}+2)p^{\delta_{0}},3(\eta_{1}+2)p^{\delta_{1}}\rbrace
		& \text{if } p^{\varsigma }-p^{\varsigma -\delta_{0}}+\eta_{0}p^{\varsigma -\delta_{0}-1}+1
		\leq \iota \leq p^{\varsigma }-p^{\varsigma -\delta_{0}}+(\eta_{0}+1)p^{\varsigma -\delta_{0}-1}\\
		& p^{\varsigma }-p^{\varsigma -\delta_{1}}+\eta_{1}p^{\varsigma -\delta_{1}-1}+1
		\leq\mu \leq p^{\varsigma }-p^{\varsigma -\delta_{1}}+(\eta_{1}+1)p^{\varsigma -\delta_{1}-1},
		\\
		
		3(\eta_{1}+2)p^{\delta_{1}}& \text{if} ~~ \iota=p^{\varsigma },\\
		&  p^{\varsigma }-p^{\varsigma -\delta_{1}}+\eta_{1}p^{\varsigma -\delta_{1}-1}+1\leq\mu \leq p^{\varsigma }-p^{\varsigma -\delta_{1}}+(\eta_{1}+1)p^{\varsigma -\delta_{1}-1}\\
		
		0,& \text{if}~~~ \iota=\mu= p^{\varsigma }.
		\end{cases}
		\end{equation*}}
\end{thm}

 To successfully develop practical quantum computers, the implementation of quantum error-correcting codes is crucial. For many years, safeguarding information from the disruptive effects of quantum noise has been a formidable challenge for researchers and engineers in the field. However, the landscape has changed significantly with remarkable breakthroughs in quantum error correction (QEC) codes. A pivotal moment in this journey was marked by the introduction of the first quantum error-correcting codes, pioneered by Peter Shor and Calderbank $\cite{calderbank1996good}$. Their contributions have paved the way for advances in protecting quantum information and improving the reliability of quantum computing systems.
 
Let $ H_{q}(C)=H_{q}(C)\otimes \cdot\cdot\cdot\otimes H_{q}(C)$ ($n$ times) be a $ q-$ dimensional Hilbert vector space (where $ q=p^{m}$).
    Then $ H_{q}^{n}(C)$ is said to be a $ q^{n}-$   dimensional Hilbert space. The definition of QEC codes is given below:
\begin{defn} \cite{rains1998quantum}
	A quantum code of length $ n$ and dimension $k$  over $ \mathbb {F}_{q}$ is defined to be a $ q^{k}$ dimensional subspace of $ H_{q}^{n}(C)$ and simply denoted by $ [[n,k,d]]_{q}$, where $ d$ is the distance of the quantum code.
\end{defn}
\begin{defn} \cite{roman1992coding}
The Singleton bound: $ \vert C \vert \leq p^{m(n-d(C)+1)}.$ If $ \vert C\vert =p^{m(n-d_{H}(C))},$ then $ C$ is said to be an AMDS code.
\end{defn}

\section{AMDS Constacyclic Codes of Length $ 4p^{\varsigma }$ over $ \mathbb{F}_{p^{m}}$}

In this section, we delve into the construction of AMDS constacyclic codes. Specifically, we utilize repeated root constacyclic codes that have a length of $4p^{\varsigma}$ over the finite field $\mathbb{F}_{p^{m}}$. This approach allows us to explore the properties and applications of these codes in greater detail.

\begin{thm}\label{thm:3.1} Assume that $ p$ is an odd prime. Let the integers $\iota, \jmath,\mu , \ell $ be such that $0\leq \iota , \jmath,\mu  ,\ell\leq p^{\varsigma}. $ 
	\begin{itemize}
		\item[(I)]  Let $ 0\leq \eta_{0},\eta_{1},\eta_{2},\eta_{3}\leq p-2$ and $ 0\leq \delta_{3}\leq \delta_{2}\leq \delta_{1}\leq \delta_{0}\leq \varsigma -1 . $ 
		Let $ C_{\iota, \jmath ,\mu,\ell }=\langle(x-1)^{\iota}(x+1)^{\jmath}(x-\delta)^{\mu}(x+\delta)^{\ell}\rangle$ be a constacyclic code of length $ 4p^{\varsigma} $ over $\mathbb{F}_{{p}^{m}}$.
		For $ 0 \leq \ell  \leq \jmath\leq\mu \leq \iota \leq p^{\varsigma },~$ and $ \delta \in \mathbb{F}_{p^m}$  such that $ \delta ^2 = -1. $ Then,  $C_{\iota, \jmath ,\mu,\ell }  $ be AMDS  if 
		\begin{itemize}
			\item[(a)] $ \iota=2, and~ \jmath =\mu=\ell =0. $ In this case $ d(C_{2,0,0,0})=2$
			\item[(b)]  $ \varsigma=1 , \jmath =\ell =0, \iota=2 ~ and ~ \mu=1.$  In this case $ d(C_{2,0,1,0})=3$ 
		\end{itemize}
		\item[(II)]
		Let $ 0\leq \eta_{0},\eta_{1},\eta_{2},\eta_{3}\leq p-2$ and $ 0\leq \delta_{3}\leq \delta_{2}\leq \delta_{1}\leq \delta_{0}\leq \varsigma -1 . $ 
		Let $ C_{\iota, \jmath ,\mu,\ell }=\langle(x-1)^{\iota}(x+1)^{\jmath}(x-\delta)^{\mu}(x+\delta)^{\ell}\rangle$ be a constacyclic code of length $ 4p^{\varsigma} $ over $\mathbb{F}_{{p}^{m}}$. 
		For $ 0 <\ell  \leq \mu \leq \jmath\leq \iota \leq p^{\varsigma },~$ and $ \delta \in \mathbb{F}_{p^m}$  such that $ \delta ^2 = -1. $ Then  $C_{\iota, \jmath ,\mu,\ell }  $ be AMDS  if
		\begin{itemize}
			\item[(a)] $ \iota=2, \jmath =\mu=\ell =0. $ In this case $ d(C_{2,0,0,0})=2$
			\item[(b)]  $  \iota = \jmath=1,  \mu=\ell =0.$  In this case $ d(C_{1,1,0,0})=2$

		\end{itemize}
		\item [(III)]
		Let $ 0\leq \eta_{0},\eta_{1},\eta_{2}\leq p-2 $ and $0\leq \delta_{2}\leq \delta_{1}\leq \delta_{0}\leq \varsigma -1 .  $ Let  $ C_{\iota, \jmath ,\mu}=\langle(x-1)^{\iota}(x+1)^{\jmath}(x^2+1)^{\mu}\rangle$ be a code of length $ 4p^{\varsigma} $ over $\mathbb{F}_{{p}^{m}}.$ 
		If $ 0\leq\mu \leq \jmath\leq \iota \leq p^{\varsigma} $ then $ C_{\iota, \jmath ,\mu}$ is an  AMDS constacyclic code if  \begin{itemize}
			\item [(a)] $ \mu=0, \jmath =1, \iota=1  $  in this case  $d(C_{1, 1, 0})=2$
			\item[(b)]  $\mu=\jmath =0, \iota=2 $ in this case $ d(C_{2,0,0})=2$
		\end{itemize} 	
		\item[(IV)] 
		Let $ 0\leq \eta_{0},\eta_{1},\eta_{2}\leq p-2 $ and $0\leq \delta_{2}\leq \delta_{1}\leq \delta_{0}\leq \varsigma -1 .  $ Let  $ C_{\iota, \jmath ,\mu}=\langle(x-1)^{\iota}(x+1)^{\jmath}(x^2+1)^{\mu}\rangle$ be a code of length $ 4p^{\varsigma} $ over $\mathbb{F}_{{p}^{m}}.$ 
		If $ 0\leq \jmath\leq\mu \leq \iota \leq p^{\varsigma} $ then $ C_{\iota, \jmath ,\mu}$ is an  AMDS constacyclic code if  \begin{itemize}
			\item [(a)] $ \mu=0, \jmath =0, \iota=2  $  in this case  $d(C_{2, 0, 0})=2$

		\end{itemize}		
		\item [(V)] 
		Let $ 0\leq {\eta}_{0},\eta_{1}\leq {p-2,}$ and $0\leq{\delta}_{0},{\delta}_{1}\leq{\varsigma-1}.$ If $C_{\iota, \jmath }=\langle(x^2-\phi)^{\iota}(x^2+\phi)^{\jmath}\rangle$ is a constacyclic code of length $ 4p^{\varsigma} $ over $\mathbb{F}_{{p}^{m}}$ then  $ C_{\iota, \jmath }$ is an AMDS constacyclic code if 
		\begin{itemize}
			\item[(a)] For $ 0\leq \jmath\leq \iota \leq p^{\varsigma}$ when $\jmath =0,\iota=1$ in this case  $d (C_{1,0})=2$  
			\item[(b)] For  $ 0\leq \iota \leq \jmath \leq p^{\varsigma},$ when $ \iota=0$ and $\jmath =1$ in this case
			$ d(C_{0,1})=2.$
		\end{itemize}   
		\item[(VI)]
		Assume that constacyclic code of length $ 4p^{\varsigma} $ over $\mathbb{F}_{{p}^{m}}$ have the form  $ C_{\iota, \mu }=\langle(x^{2}+\gamma \zeta^{\frac{2p^{m-r}+3p^m-1}{8}}x+ \zeta^{\frac{2p^{m-r}+p^{m}-1}{4}})^{\iota}(x^{2}-\gamma \zeta^{\frac{2p^{m-r}+3p^m-1}{8}}x+ \zeta^{\frac{2p^{m-r}+p^{m}-1}{4}})^{\mu})  \rangle$,  and $  \gamma \in \mathbb {F}_{p^{m}}$ such that $ \gamma ^{2}=-2.$ Then $C_{\iota,\mu} $ is an AMDS constacyclic code if 
		\begin{itemize}	
	\item[(a)]	For    $ 0 \leq\mu \leq \iota \leq p^{\varsigma },$  when $ \iota=1, \mu=0$ in this case $ d(C_{1,0})=2.$
	\item[(b)] For $ 0 \leq \iota \leq\mu \leq p^{\varsigma },$  when $ \iota=0, \mu=1$ in this case $ d(C_{0,1})=2.$
		\end{itemize}	
		
		\item[(VII)]Let $ p$ be an odd prime, $m$ be a positive integer, $ \psi $ be nonzero element of $ \mathbb{F}_{p^{m}}$.
		Assume that $ 0 \leq \eta_{0}\leq p-2$ and $ 0 \leq \delta_{0}\leq \varsigma -1 $ and $ C_{\iota}=\langle(x^4+\psi)^{\iota}\rangle$ is a constacyclic code of length $ 4 p^{\varsigma }$ over $ \mathbb{F}_{p^{m}},$ where $ 0\leq \iota \leq p^{\varsigma }$ and $ x^4+ \psi$ is irreducible over $ \mathbb{F}_{p^{m}}$. Then $ C_{\iota}$ is not an AMDS code of length $ 4p^{\varsigma  }$ over $ \mathbb{F}_{p^{m}}$
	\end{itemize}
\end{thm}  
\begin{pro}: Here, we give proof of part (I).\\

	\textbf{Case 1:} Let  $  \iota = \jmath=\mu=\ell =0 .$ Using Theorem \ref{thm:four} , we have $ d(C_{\iota, \jmath ,\mu,\ell })=1.$ In this case, we note that $ \iota + \jmath +\mu+\ell =0\neq d(C_{0,0,0,0}).$ Therefore, the constacyclic code $ C_{0,0,0,0}$ is not an AMDS Code of length $ 4p^{\varsigma }.$\\
	\textbf{Case 2:} Let $ \jmath =\mu=\ell =0$ and $ 0 < \iota \leq p^{\varsigma }$ or $ \ell =0$ and $ 0 < \jmath \leq \mu \leq \iota \leq p^{\varsigma -1}.$
	Using Theorem \ref{thm:four}, we see that $ d(C_{\iota, \jmath ,\mu,\ell })=2.$ In this case, we note  that if $ \iota=2, \jmath =\mu=\ell =0,$ then $ \iota + \jmath +\mu+\ell =2=d(C_{\iota ,\jmath,k,\ell }).$ Therefore, the constacyclic code $ C_{2,0,0,0}$is an AMDS code  of length  $4p^{\varsigma}$  over $ \mathbb {F}_{p^m}.$\\
	\textbf{Case 3:} Let $ \jmath =\ell =0, 0<\mu\leq p^{\varsigma}$ and $ p^{\varsigma -1}<\iota \leq p^{\varsigma},$ or $ \ell =0, 0<\jmath \leq k\leq\mu \leq 2p^{\varsigma -1}$
	and $ p^{\varsigma -1} < \iota \leq 2p^{\varsigma -1}.$ Using Theorem \ref{thm:four}, we have $ d(C_{\iota, \jmath ,\mu,\ell })=3.$ Therefore,
	$ \iota + \jmath +\mu+\ell\geq   p^{\varsigma -1}+2.$ For $ \varsigma=1, \iota=2, \mu=1, \jmath =\ell =0$  we get $d(C_{2,0,1,0})=3.$ Hence AMDS. But for $ s>1, C_{\iota, \jmath ,\mu,\ell }$ is not AMDS code.\\
	\textbf{Case 4:}Let $ \ell =0, 0< \jmath \leq \mu \leq p^{\varsigma}$ and $ 2p ^{\varsigma -1}<\iota \leq p^{\varsigma }.$ Using Theorem \ref{thm:four} we see
	that $ d(C_{\iota, \jmath ,\mu,\ell })=4$ In this case,we note that $ \iota + \jmath +\mu+\ell\geq p^{\varsigma -1}+ 3.$ Which is not possible. Hence, the constacyclic code $ C_{\iota, \jmath,\mu,\ell } $ is not an AMDS code.\\
	\textbf{Case 5:} Let $ p^{\varsigma }-p^{\varsigma -\delta_{0}}+\eta_{0}p^{\varsigma -\delta_{0}-1}+1
	\leq \iota \leq p^{\varsigma }-p^{\varsigma -\delta_{0}}+(\eta_{0}+1)p^{\varsigma -\delta_{0}-1},~~~
	p^{\varsigma }-p^{\varsigma -\delta_{1}}+\eta_{1}p^{\varsigma -\delta_{1}-1}+1
	\leq \jmath\leq p^{\varsigma }-p^{\varsigma -\delta_{1}}+(\eta_{1}+1)p^{\varsigma -\delta_{1}-1},~~
	p^{\varsigma }-p^{\varsigma -\delta_{2}}+\eta_{2}p^{\varsigma -\delta_{2}-1}+1
	\leq k \leq p^{\varsigma }-p^{\varsigma -\delta_{2}}+(\eta_{2}+1)p^{\varsigma -\delta_{2}-1} ~~and~~
	p^{\varsigma }-p^{\varsigma -\delta_{3}}+\eta_{3}p^{\varsigma -\delta_{3}-1}+1
	\leq \ell  \leq p^{\varsigma }-p^{\varsigma -\delta_{3}}+(\eta_{3}+1)p^{\varsigma -\delta_{3}-1}$ with
	$ 0\leq \eta_{0},\eta_{1},\eta_{2},\eta_{3}\leq p-2 $ and $0\leq\delta_{3}\leq \delta_{2}\leq \delta_{1}\leq \delta_{0}\leq \varsigma -1 .$ Using theorem \ref{thm:four} we have
	$  d(C_{\iota, \jmath ,\mu,\ell })= min \lbrace(\eta_{0}+2)p^{\delta_{0}},2(\eta_{1}+2)p^{\delta_{1}}, 3(\eta_{2}+2)p^{\delta_{2}},4(\eta_{3}+2)p^{\delta_{3}}\rbrace $ This implies that
	\begin{align*}         
	\iota + \jmath +k+\ell &\geq 4p^{\varsigma}-p^{\varsigma -\delta_{0}}-p^{\varsigma -\delta_{1}}-p^{\varsigma -\delta_{2}}-p^{\varsigma -\delta_{3}}+\eta_{0}p^{\varsigma -\delta_{0}-1}+
	\eta_{1}p^{\varsigma -\delta_{1}-1}+\eta_{2}p^{\varsigma -\delta_{2}-1}+\eta_{3}p^{\varsigma -\delta_{3}-1}+4\\ 
	&= 4p^{\varsigma }-4p^{\varsigma -\delta_{0}}+4\eta_{0}p^{\varsigma -\delta_{0}-1}+4
	(equality~ holds~ when~ \eta_{0}=\eta_{1} =\eta_{2}=\eta_{3}~ and~ \delta_{0}=\delta_{1}=\delta_{2}=\delta_{3})\\
	&=4(p^{\varsigma -\delta_{0}}(p^{\delta_{0}}-1))+\eta_{0}p^{\varsigma -\delta_{0}-1}+1\\
	&\geq 4(p(p^{\delta_{0}}-1)+\eta_{0}+1)(equality~ holds~ when~ s-\delta_{0}=1)\\
	&\geq 4((\eta_{0}+2)(p^{\delta_{0}}-1)+\eta_{0}+1)(equality~ holds~ when~ p=\eta_{0}+2)\\
	&=4(\eta_{0}+2)p^{\delta_{0}}-4\\
	&=(\eta_{0}+2)p^{\delta_{0}}+3(\eta_{0}+2)p^{\delta_{0}}-4\\
	&>(\eta_{0}+2)p^{\delta_{0}}~~~~~~~(as ~~ 3(\eta_{0}+2)p^{\delta_{0}}-4 ) > 0) \\
	& >min\lbrace(\eta_{0}+2)p^{\delta_{0}},2(\eta_{1}+2)p^{\delta_{1}},  3(\eta_ 
	{2}+2)p^{\delta_{2}}, 4(\eta_{3}+2)p^{\delta_{3}}\rbrace\\
	&=d(C_{\iota, \jmath ,\mu,\ell })       
	\end{align*}
	Therefore, the constacyclic code $ C_{\iota, \jmath ,\mu,\ell }$ is not an AMDS code.\\
	\textbf{Case 6: } Let $\iota=p^{\varsigma },p^{\varsigma }-p^{\varsigma -\delta_{1}}+\eta_{1}p^{\varsigma -\delta_{1}-1}+1\leq \jmath\leq p^{\varsigma }-p^{\varsigma -\delta_{1}}+(\eta_{1}+1)p^{\varsigma -\delta_{1}-1},\\
	p^{\varsigma }-p^{\varsigma -\delta_{2}}+\eta_{2}p^{\varsigma -\delta_{2}-1}+1\leq k\leq p^{\varsigma }-p^{\varsigma -\delta_{2}}+(\eta_{2}+1)p^{\varsigma -\delta_{2}-1},
	p^{\varsigma }-p^{\varsigma -\delta_{3}}+\eta_{3}p^{\varsigma -\delta_{3}-1}+1
	\leq \ell  \leq p^{\varsigma }-p^{\varsigma -\delta_{3}}+(\eta_{3}+1)p^{\varsigma -\delta_{3}-1} $ with 
	$ 0 \leq \eta_{1},\eta_{2},\eta_{3}\leq p-2$ and $ 0 \leq \delta_{3}\leq \delta_{2}\leq \delta_{1}\leq \varsigma -1 .$ Using Theorem \ref{thm:four} we have $ d(C_{\iota, \jmath, \mu  ,\ell })= min\lbrace 2(\eta_{1}+2)p^{\delta_{1}},  3(\eta_{2}+2)p^{\delta_{2}}, 4(\eta_{3}+2)p^{\delta_{3}}\rbrace. $ Since
	\begin{align*}        
	\iota + \jmath +\mu+\ell&\geq 4p^{\varsigma}-p^{\varsigma -\delta_{1}}-p^{\varsigma -\delta_{2}}-p^{\varsigma -\delta_{3}}+\eta_{1}p^{\varsigma -\delta_{1}-1}+\eta_{2}p^{\varsigma -\delta_{2}-1}+\eta_{3}p^{\varsigma -\delta_{3}-1}+3 \\
	&=4p^{\varsigma }-3p^{\varsigma -\delta_{1}}+3\eta_{1}p^{\varsigma -\delta_{1}-1}+3(equality~holds~when~\eta_{1}=\eta_{2}=\eta_{3}~and~\delta_{1}=\delta_{2}=\delta_{3})\\
	&=p^{\varsigma }+3((p^{\varsigma -\delta_{1}}(p^{\delta_{1}}-1))+\eta_{1}p^{\varsigma -\delta_{1}-1}+1)\\
	&\geq p^{\varsigma }+3(p(p^{\delta_{1}}-1)+\eta_{1}+1)(equality~holds~when~s-\delta_{1}=1)\\
	&\geq p^{\varsigma }+3((\eta_{1}+2)(p^{\delta_{1}}-1)+\eta_{1}+1)(equality~holds~when~p=\eta_{1}+2)\\
	&=p^{\varsigma }+2((\eta_{1}+2)p^{\delta_{1}} +(\eta_{1}+2)p^{\delta_{1}} -3\\
	& >2(\eta_{1}+2)p^{\delta_{1}},   ( since ~~ (\eta_{1}+2)p^{\delta_{1}} -3>0) \\
	&>min\lbrace 2((\eta_{1}+2)p^{\delta_{1}}),  3(\eta_{2}+2)p^{\delta_{2}},4(\eta_{3}+2)p^{\delta_{3}}  \rbrace\\
	&\neq d(C_{\iota, \jmath ,\mu,\ell })                   
	\end{align*}
	Therefore,the constacyclic code $ C_{\iota, \jmath ,\mu,\ell }$ is not an AMDS constacyclic code.\\
	\textbf{Case 7:} Let $  \iota = \jmath=p^{\varsigma },  p^{\varsigma }-p^{\varsigma -\delta_{2}}+\eta_{2}p^{\varsigma -\delta_{2}-1}+1\leq k\leq p^{\varsigma }-p^{\varsigma -\delta_{2}}+(\eta_{2}+1)p^{\varsigma -\delta_{2}-1},
	p^{\varsigma }-p^{\varsigma -\delta_{3}}+\eta_{3}p^{\varsigma -\delta_{3}-1}+1
	\leq \ell  \leq p^{\varsigma }-p^{\varsigma -\delta_{3}}+(\eta_{3}+1)p^{\varsigma -\delta_{3}-1}, $  with $ 0\leq\eta_{2},\eta_{3}\leq p-2$ and $ 0\leq \delta_{3}\leq \delta_{2}\leq \varsigma -1 .$ Using Theorem \ref{thm:four} we see that $d(C_{\iota, \jmath ,\mu,\ell })=min\lbrace 3(\eta_{2}+2)p^{\delta_{2}}, 4(\eta_{3}+2)p^{\delta_{3}}     \rbrace $. Consider
	\begin{align*}
	\iota + \jmath +k+\ell &\geq 4p^{\varsigma}-p^{\varsigma -\delta_{2}}-p^{\varsigma -\delta_{3}}+\eta_{2}p^{\varsigma -\delta_{2}-1}+\eta_{3}p^{\varsigma -\delta_{3}-1}+2 ,\\
	&= 4p^{\varsigma }-2p^{\varsigma -\delta_{2}}+2\eta_{2}p^{\varsigma -\delta_{2}-1}+2 (equality~holds~when~\eta_{2}=\eta_{3}~~ and~~ \delta_{2}=\delta_{3})\\
	&=2p^{\varsigma -\delta_{2}}p^{\delta_{2}}+2(p^{\varsigma -\delta_{2}}(p^{\delta_{2}}-1)+\eta_{2}p^{\varsigma -\delta_{2}-1}+1)\\
	&\geq2p(p^{\delta_{2}})+2(p(p^{\delta_{2}}-1)+\eta_{2}+1)(equality~~holds~~when~s-\delta_{2}=1 )\\
	&\geq 2(\eta_{2}+2)(p^{\delta_{2}})+2[(\eta_{2}+2)(p^{\delta_{2}}-1)+\eta_{2}+1](equality~~holds~~when~p=\eta_{2}+2 )\\
	&= 2(\eta_{2}+2)(p^{\delta_{2}})+2[(\eta_{2}+2)p^{\delta_{2}}-1]\\
	&= 3(\eta_{2}+2)(p^{\delta_{2}})+(\eta_{2}+2)(p^{\delta_{2}})-2\\ 
	&>3(\eta_{2}+2)(p^{\delta_{2}})\\
	&>min\lbrace 3(\eta_{2}+2)(p^{\delta_{2}}), 4(\eta_{2}+2)(p^{\delta_{3}})  \rbrace\\
	&=d(C_{\iota, \jmath ,\mu,\ell })
	\end{align*}
	Therefore, the constacyclic code $ C_{\iota, \jmath ,\mu,\ell }$ is not an AMDS code.\\
	\textbf{Case 7:} Let $  \iota = \jmath=\mu= p^{\varsigma },    p^{\varsigma }-p^{\varsigma -\delta_{3}}+\eta_{3}p^{\varsigma -\delta_{3}-1}+1
	\leq \ell  \leq p^{\varsigma }-p^{\varsigma -\delta_{3}}+(\eta_{3}+1)p^{\varsigma -\delta_{3}-1}$ with $ 0\leq \eta_{3}\leq p-2$ and $0\leq \delta_{3}\leq \varsigma -1 .$ Using Theorem \ref{thm:four} we see that $ d({C_{\iota, \jmath ,\mu,\ell }})= 4(\eta_{2}+2)(p^{\delta_{3}}). $ Consider  
	\begin{align*}
	\iota + \jmath +\mu+\ell&\geq 4p^{\varsigma }-p^{\varsigma -\delta_{3}}+\eta_{3}p^{\varsigma -\delta_{3}-1}+1\\
	&=3p^{\varsigma }+p^{\varsigma -\delta_{3}}(p^{\delta_{3}}-1)+\eta_{3}p^{\varsigma -\delta_{3}-1}+1\\
	&\geq 3p^{\eta_{3}+1}+p(p^{\delta_{3}}-1)+\eta_{3}+1 (equality~~holds~~when~s-\delta_{3}=1 )\\
	&\geq 3(\eta_{3}+2)p^{\eta_{3}}+(\eta_{3}+2)(p^{\delta_{3}}-1)+\eta_{3}+1 (equality~~holds~~when~p=\eta_{3}+2 )\\
	&=4 (\eta_{3}+2)p^{\delta_{3}}-1\\
	&=d(C_{\iota, \jmath ,\mu,\ell })-1        
	\end{align*}
	\textbf{Case 8:} Let $  \iota = \jmath=\mu=\ell =p^{\varsigma }.$ Using Theorem \ref{thm:four}, we see that $ d(C_{\iota, \jmath ,k,\ell })=0.$ In this case, $ \iota + \jmath +\mu+\ell =4p^{\varsigma }\neq d(C_{\iota, \jmath ,\mu,\ell }.)$ Hence, it is not an AMDS code.\\
    Working similarly to part (I) and using Theorems \ref{2.2}, \ref{thm:one}, \ref{thm:two}, \ref{thm:three}, and \ref{thm:five}, we obtain results for the other parts as well.
\end{pro}
\begin{remark} Following part (I) of  Theorem \ref{thm:3.1}, we can determine all AMDS constacyclic codes of length $ 4p^{\varsigma}$ over $ \mathbb{F}_{p^m}$ with 	$\mu \leq \ell  \leq \jmath \leq \iota, \ell  \leq \mu \leq \iota \leq \jmath,  \mu \leq \ell  \leq \iota \leq \jmath ,  \jmath \leq \iota \leq \ell  \leq \mu, \iota \leq \jmath \leq \ell \leq \mu , \iota \leq \jmath\leq \mu \leq \ell ,~~ and~~~ \jmath \leq \iota \leq \mu \leq \ell .$
	\end{remark}

\begin{example}
	Let $ p=5, \varsigma=1,m =1.$ Then we have $ x^{20}-1=(x-1)^5(x+1)^5(x-2)^5(x+2)^5.$\\ Let $ C_{2,0,0,0}= \langle(x-1)^2\rangle$  is a code with  parameters [20,18,2] .\\
	Let $ C_{2,0,1,0}=  \langle(x-1)^2 (x-2)\rangle$ is a code with  parameters [20,17,3].\\
Let $ C_{1,1,0,0}= \langle(x-1)(x+1)\rangle$  is a code with  parameters [20,18,2].
	  Using Theorem \ref{thm:3.1}  part (I) and part(II), $ C_{2,0,0,0}$  $ C_{2,0,1,0}$ and $ C_{1,1,0,0}$ are AMDS constacyclic codes.	
\end{example}

\begin{example}
	Let $ p=7,\varsigma=1,m=1.$ Then we have $ x^{28}-1=(x-1)^7(x+1)^7(x^2+1)^7.$\\ Let $ C_{1,1,0}= \langle(x-1)(x+1)\rangle$  is a code with  parameters [28,26,2]. \\

	Let $C_{2,0,0}=\langle (x-1)^2\rangle $	 is a code with parameters [28,26,2].\\

	Using Theorem \ref{thm:3.1} part (III), $ C_{1,1,0}$ and $C_{2,0,0}$ are AMDS constacyclic code.
\end{example}
\begin{example}
	Let $ p=5,~~~\varsigma=1,~~~m=1.$ Then we have $ x^{20}-2=(x^2+2)^5(x^2+3)^5.$ Let $ C_{0,1}= \langle(x^2-2)\rangle$ is a code with  parameters [20,18,2] .Using Theorem \ref{thm:3.1} part (V), $ C_{0,1}$ is an AMDS constacyclic code.	
\end{example}

\begin{example}
	Let $ p=3,~~~\varsigma=1,~~~m=1.$ Then we have $ x^{12}-2=(x^2+x+2)^3(x^2+2x+2)^3.$ Let $ C_{1,0}= \langle(x^2+x+2)\rangle$  is a code with parameters [12,10,2] .Using Theorem \ref{thm:3.1} part(VI) , $ C_{1,0}$ is an AMDS constacyclic code.
\end{example}

\section{Quantum AMDS Constacyclic  Codes from CSS constructions}

In this section, we provide a comprehensive overview of the essential concepts of quantum codes, including the CSS construction method used to define them. Following this introduction, we will delve into the necessary and sufficient conditions required for the existence of dual-containing constacyclic codes of a specific length $4p^{\varsigma} $ over the finite field $\mathbb{F}_{p^m} $. 
Building upon these dual-containing repeated root constacyclic codes, we will employ the CSS construction technique to systematically generate and identify all possible quantum codes of the specified length $4p^{\varsigma}$  over the field  $\mathbb{F}_{p^m} $. This process not only enhances our understanding of quantum error correction but also expands the repertoire of codes available for practical applications in quantum computing.

\begin{thm}\label{4.1} \cite{calderbank1996good}
	(CSS Construction)Assume that $ C_{1}$ is a linear code which has parameters $ [n,k_{1},d_{1}]_{q}$ and $ C_{2}$ is a linear code over $ \mathbb{F}_{p^m}$ which has parameters $ [n,k_{2},d_{2}]_{q}$ satisfying $ C_{2}\subseteq C_{1}.$ Then there exists a QEC code with the parameters $ [[n, k_{1}-k_{2}, min{d_{1},d_{2}^{\perp}}]]_{q}$ ($d_{2}^{\perp} $ is the Hamming distance of the dual code $ C_{2}^{\perp}$). Specially, if $C_{2}=C_{1}^{\perp}, $ then there exists a QEC code having the parameters $ [[n,2k_{1}-n,d_{1}]]_{q}.$
\end{thm} 
\begin{thm}\label{4.2}
	(Quantum Singleton Bound)\cite{grassl2004optimal} Let $ C=[[n,k,d]]_{q}$ be a QEC code. Then $ 2d+k\leq n+2.$  
\end{thm}
If $ 2d=n-k$, then $C$ is called a qAMDS code.
\begin{thm}\label{thm:4.1}
	Let $ p$ be an odd prime and integers $ \iota, \jmath  ,\mu,\ell $  run over the set $ \{0,1,...,p^{\varsigma}\}$,
	\begin{itemize}
		\item [I)] 
		Let $ 0 \leq \ell  \leq\mu \leq \jmath\leq \iota \leq p^{\varsigma },  $ and $ C_{\iota, \jmath ,\mu,\ell }=\langle(x-1)^{\iota}(x+1)^{\jmath}(x-\delta)^{\mu}(x+\delta)^{\ell}\rangle$ be a constacyclic code of length $ 4p^{\varsigma} $ over $\mathbb{F}_{{p}^{m}} $ with dual  $ C^{\perp}_{\iota, \jmath ,\mu,\ell }=\langle(x-1)^{p^{\varsigma}- \iota}(x+1)^{p^{\varsigma}- \jmath}(x-\delta)^{p^{\varsigma}-\mu}(x+\delta)^
		{p^{\varsigma}-l}\rangle .$ 
		
		 If $   \iota = \jmath=1, \mu=\ell =0$ then there exists a qAMDS code with parameters $[[4p^{\varsigma}, 4p^{\varsigma}-4,2]]_{p^m}.$ If $   \iota=2, \jmath =0, \mu=\ell =0$ then there exists a qAMDS code with parameters $   [[4p^{\varsigma}, 4p^{\varsigma}-4,2]]_{p^m}     .$  
		\item[II)] Let $ 0 <\ell  \leq  \jmath\leq\mu \leq \iota \leq p^{\varsigma },  $ and $ C_{\iota, \jmath ,\mu,\ell }=\langle(x-1)^{\iota}(x+1)^{\jmath}(x-\delta)^{\mu}(x+\delta)^{\ell}\rangle$ be a constacyclic code of length $ 4p^{\varsigma} $ over $\mathbb{F}_{{p}^{m}}$ with dual as $ C^{\perp}_{\iota, \jmath ,\mu,\ell }=\langle(x-1)^{p^{\varsigma}- \iota}(x+1)^{p^{\varsigma}- \jmath}(x-\delta)^{p^{\varsigma}-\mu}(x+\delta)^{p^{\varsigma}-l}\rangle .$ If $  \varsigma=1,\iota=2, \jmath =\mu=\ell =0$ then there exists a qAMDS code with parameters $[[4p^{\varsigma}, 4p^{\varsigma}-4,2]]_{p^m}.$ If $  \varsigma=1, \iota=2, \jmath =\ell =0, \mu=1$ then there exists a qAMDS code with parameters $[[4p^{\varsigma}, 4p^{\varsigma}-6,3]]_{p^m}.$
		\item [III)]  Let $ 0\leq\mu \leq \jmath\leq \iota \leq p^{\varsigma}$ and $ C_{\iota, \jmath ,\mu}=\langle(x-1)^{\iota}(x+1)^{\jmath}(x^2+1)^{\mu}\rangle$ be a code of length $ 4p^{\varsigma} $ over $\mathbb{F}_{{p}^{m}}$ with dual as $ C^{\perp}_{\iota, \jmath ,\mu}=\langle(x-1)^{p^{\varsigma}- \iota}(x+1)^{p^{\varsigma}- \jmath}(x^2+1)^{p^{\varsigma}-\mu}\rangle.$\\	
		If $  \iota = \jmath=1, \mu=0 $ then there exists a quantum AMDS code with parameters $ [[4p^{\varsigma}, 4p^{\varsigma}-4,2]]_{p^m}.$
		If $ \iota=2, \jmath =0, \mu=0$ then there exists a quantum AMDS code with parameters $ [[4p^{\varsigma}, 4p^{\varsigma}-4,2]]_{p^m}$
		\item[IV)] Let $ 0\leq \jmath\leq\mu \leq \iota \leq p^{\varsigma}$ and $ C_{\iota, \jmath ,\mu}=\langle(x-1)^{\iota}(x+1)^{\jmath}(x^2+1)^{\mu}\rangle$ be a code of length $ 4p^{\varsigma} $ over $\mathbb{F}_{{p}^{m}}$ with dual as $ C^{\perp}_{\iota, \jmath ,\mu}=\langle(x-1)^{p^{\varsigma}- \iota}(x+1)^{p^{\varsigma}- \jmath}(x^2+1)^{p^{\varsigma}-\mu}\rangle.$ If $ \iota=2,\jmath =0 ~and ~\mu=0,$ then there exists a qAMDS code with parameters $[[4p^{\varsigma}, 4p^{\varsigma}-4,2]]_{p^m}.$ If $ \varsigma=1,\iota=1, \jmath =0 ~~and ~~ \mu=1$ then there exists a qAMDS code with parameters $[[4p^{\varsigma}, 4p^{\varsigma}-6,3]]_{p^m}.$
		\item[V)]	
		let $ 0\leq \jmath\leq \iota \leq p^{\varsigma} $ and $C_{\iota, \jmath }=\langle(x^2-\phi)^{\iota}(x^2+\phi)^{\jmath}\rangle$ is a constacyclic code of length $ 4p^{\varsigma} $ over $\mathbb{F}_{{p}^{m}}$ with dual $ C^{\perp}_{i,j}= \langle(x^2-\phi)^{p^{\varsigma}- \iota}(x^2+\phi)^{p^{\varsigma}- \jmath}\rangle $ then there exists a quantum AMDS
		code having the parameters  $ [[4p^{\varsigma},4p^{\varsigma}-4,2]]_{p^m}$.
		\item[VI)] 	 
		Assume that constacyclic code of length $ 4p^{\varsigma} $ over $\mathbb{F}_{{p}^{m}}$ have the form  $ C_{\iota, \mu }=\langle(x^{2}+\gamma \zeta^{\frac{2p^{m-r}+3p^m-1}{8}}x+ \zeta^{\frac{2p^{m-r}+p^{m}-1}{4}})^{\iota}(x^{2}-\gamma \zeta^{\frac{2p^{m-r}+3p^m-1}{8}}x+ \zeta^{\frac{2p^{m-r}+p^{m}-1}{4}})^{\mu})  \rangle$,   $  \gamma \in \mathbb {F}_{p^{m}}$ such that $ \gamma ^{2}=-2$ and with dual as $C^{\perp}_{\iota, \mu }= \langle(x^{2}+\gamma \zeta^{\frac{2p^{m-r}+3p^m-1}{8}}x+ \zeta^{\frac{2p^{m-r}+p^{m}-1}{4}})^{p^{\varsigma}- \iota}(x^{2}-\gamma \zeta^{\frac{2p^{m-r}+3p^m-1}{8}}x+ \zeta^{\frac{2p^{m-r}+p^{m}-1}{4}})^{p^{\varsigma}-\mu}) . $ 
				If    $ 0 \leq\mu \leq \iota \leq p^{\varsigma },$ then $ C_{\iota, \mu }$  then there exists a quantum  AMDS constacyclic code when $ \iota=1, \mu=0 $ with parameters $[[4p^{\varsigma}, 4p^{\varsigma}-2, 2]]_{p^m}. $

		\end{itemize}
\end{thm}
\begin{pro}
	\begin{itemize}
		
		\item [I)] Let $ 0 \leq \ell  \leq\mu \leq \jmath\leq \iota \leq p^{\varsigma },  $ and $ C_{\iota, \jmath ,\mu,\ell}=\langle(x-1)^{\iota}(x+1)^{\jmath}(x-\delta)^{\mu}(x+\delta)^{\ell}\rangle$ be a constacyclic code of length $ 4p^{\varsigma} $ over $\mathbb{F}_{{p}^{m}}$ with dual as $ C^{\perp}_{\iota, \jmath ,\mu,\ell}=\langle(x-1)^{p^{\varsigma}- \iota}(x+1)^{p^{\varsigma}- \jmath}(x-\delta)^{p^{\varsigma}-\mu}(x+\delta)^{p^{\varsigma}-l}\rangle .$ Now, if $ C^{\perp}_{\iota, \jmath ,\mu,\ell} \subseteq C_{\iota, \jmath ,\mu,\ell},$ then 
		$p^{\varsigma}- \iota \geq \iota,~ p^{\varsigma}- \jmath \geq \jmath,~ p^{\varsigma}- \mu  \geq  \mu  $ and $ p^{\varsigma}-l \geq l .$ This means $ 0 \leq \iota , \jmath, \mu ,\ell  \leq  \frac{p^{\varsigma}}{2}.$ Let		
		 $ C_{\iota, \jmath ,\mu,\ell } $ be an AMDS constacyclic code of length $ 4p^{\varsigma}$ over $\mathbb{F}_{{p}^{m}}$ with parmeters $ [4p^{\varsigma}, k_{\iota, \jmath ,\mu,\ell },d(C_{\iota,\jmath,\mu,\ell })]_{p^m}$ which satisfies the condition that
		$C^{\perp}_{\iota, \jmath ,\mu,\ell }\subseteq C_{\iota, \jmath ,\mu,\ell }  .$ Thus $ k_{\iota, \jmath ,\mu,\ell }=4p^{\varsigma}-d(C_{\iota, \jmath ,\mu,\ell })$ and $ 0 \leq \iota , \jmath, \mu ,\ell  \leq  \frac{p^{\varsigma}}{2}.$  Now by using CSS construction as in Theorem \ref{4.1} there exists a QEC code, say $ M_{\iota,\jmath,\mu,\ell }$ having parameters $ [[4p^{\varsigma}, 2k_{\iota, \jmath ,\mu,\ell }-4p^{\varsigma}, d(C_{\iota, \jmath ,\mu,\ell })]]_{p^m}.$ Here the dimension of $ M_{\iota, \jmath ,\mu,\ell }$  is $2k_{\iota, \jmath ,\ell ,\mu}-4p^{\varsigma} $ and from above we have $ 2k_{\iota, \jmath ,\mu,\ell }=8p^{\varsigma}- 2d(C_{\iota,\jmath,\mu,\ell })$ Thus the dimension of this QEC code $ M_{\iota,\jmath,\mu,\ell }$ is $ 4p^{\varsigma}-2d(C_{\iota,\jmath,\mu,\ell }).$ By Theorem \ref{4.2} it is clear that QEC code $ M_{\iota, \jmath ,\mu,\ell }$  attains quantum singleton bound. Hence  $ M_{\iota, \jmath ,\mu,\ell }$ is a quantum AMDS code having  parameters  $[4p^{\varsigma}, 4p^{\varsigma}-2d(C_{\iota, \jmath ,\mu,\ell }), d(C_{\iota, \jmath ,\mu,\ell })]. $ Thus we conclude that if $ C_{\iota, \jmath ,\mu,\ell }$  is an AMDS constacyclic code of length $ 4p^{\varsigma}$ over $\mathbb{F}_{{p}^{m}}$ with parameters
		$ [4p^{\varsigma}, k_{\iota, \jmath ,\mu,\ell },d(C_{\iota, \jmath ,\mu,\ell })]_{p^m}$ and  $C^{\perp}_{\iota, \jmath ,\mu,\ell }\subseteq C_{\iota, \jmath ,\mu,\ell } $, then there exists  a quantum AMDS code $ M_{\iota, \jmath ,\mu,\ell }$ having the parameters $[4p^{\varsigma}, 4p^{\varsigma}-2d(C_{\iota, \jmath ,\mu,\ell }), d(C_{\iota, \jmath ,\mu,\ell })]. $ Thus we have the following possibilities;\\
		If $  \iota = \jmath=1, \mu=\ell =0 $ then by using Theorm \ref{thm:3.1} $ C_{1,1,0,0}$ is an AMDS constacyclic code having the parameters $ [4p^{\varsigma},4p^{\varsigma}-2 ,2].$ Also $ C^{\perp}_{1,1,0,0} \subseteq C_{1,1,0,0}.$ Thus, by the above argument, there exists a quantum AMDS code having the parameters $ [[4p^{\varsigma},4p^{\varsigma}-4,2]]_{p^m}.$\\
		If $ \iota=2,\jmath =0,\mu=\ell =0$ then by using Theorem \ref{thm:3.1}
		$ C_{2,0,0,0}$ is an AMDS constacyclic code having the parameters $ [[4p^{\varsigma}, 4p^{\varsigma}-2,2]].$ Also note that  $ C^{\perp}_{2,0,0,0} \subseteq C_{2,0,0,0}.$ Thus, by similiar argument as above, there exists a quantum AMDS code having the parameters $ [[4p^{\varsigma},4p^{\varsigma}-4,2]]_{p^m}.$\\

		\item[III)]  In this case  $ 0\leq\mu \leq \jmath\leq \iota \leq p^{\varsigma}$ and $ C_{\iota, \jmath ,\mu}=\langle(x-1)^{\iota}(x+1)^{\jmath}(x^2+1)^{\mu}\rangle$ be a code with dual as $ C^{\perp}_{\iota, \jmath ,\mu}=\langle(x-1)^{p^{\varsigma}- \iota}(x+1)^{p^{\varsigma}- \jmath}(x^2+1)^{p^{\varsigma}-\mu}\rangle.$
		Now, if $ C^{\perp}_{\iota, \jmath ,\mu} \subseteq C_{\iota, \jmath ,\mu},$ then $p^{\varsigma}- \iota \geq \iota,~ p^{\varsigma}- \jmath \geq \jmath  $ and $ p^{\varsigma}- \mu  \geq  \mu.$ This means $ 0 \leq \iota , \jmath ,  \mu  < \frac {p^{\varsigma}}{2}.$
		 Let $ C_{\iota, \jmath , \mu }$ be an AMDS constacyclic code of length $ 4p^{\varsigma}$ over $\mathbb{F}_{{p}^{m}}$
		with parmeters\\  $ [4p^{\varsigma}, k_{\iota, \jmath ,\mu},d(C_{\iota,\jmath,\mu})]_{p^m}$ which satisfies the condition that
		$C^{\perp}_{\iota, \jmath ,\mu}\subseteq C_{\iota, \jmath ,\mu}  .$ Thus $ k_{\iota, \jmath ,\mu}=4p^{\varsigma}-d(C_{\iota, \jmath ,\mu})$ and $ 0 \leq \iota , \jmath,\mu \leq \frac {p^{\varsigma}}{2}$ and $ 0\leq \jmath\leq \frac{p^{\varsigma}}{2}.$ Now by using CSS construction as in Theorem \ref{4.1} there exists a QEC code, say $ M_{\iota,\jmath,\mu}$ having parameters $ [[4p^{\varsigma}, 2k_{\iota, \jmath ,\mu}-4p^{\varsigma}, d(C_{\iota, \jmath ,\mu})]]_{p^m}.$ Here the dimension of the obtained code is $2k_{\iota, \jmath ,\mu}-4p^{\varsigma} $ and from above we have $ 2k_{\iota, \jmath ,\mu}=8p^{\varsigma}- 2d(C_{\iota,\jmath,\mu})$ Thus the dimension of this QEC code $ M_{\iota,\jmath,\mu}$ is $ 4p^{\varsigma}-2d(C_{\iota,\jmath,\mu}).$ By Theorem \ref{4.2} it is clear that QEC code $ M_{\iota, \jmath ,\mu}$  attains quantum singleton bound. Hence  $ M_{\iota, \jmath ,\mu}$ is a quantum AMDS code having 
		parameters  $[4p^{\varsigma}, 4p^{\varsigma}-2d(C_{\iota, \jmath ,\mu}), d(C_{\iota,\jmath,\mu})] .$ Thus we conclude that if $ C_{\iota, \jmath ,\mu}$  is an AMDS constacyclic code of length $ 4p^{\varsigma}$ over $\mathbb{F}_{{p}^{m}}$ with parameters
		$ [4p^{\varsigma}, k_{\iota, \jmath ,\mu},d(C_{\iota, \jmath ,\mu})]_{p^m}$ and  $C^{\perp}_{\iota, \jmath ,\mu}\subseteq C_{\iota, \jmath ,\mu} $, then there exists  a quantum AMDS code $ M_{\iota, \jmath ,\mu}$ having the parameters $[4p^{\varsigma}, 4p^{\varsigma}-2d(C_{\iota, \jmath ,\mu}), d(C_{\iota, \jmath ,\mu})]. $ Thus we have the following possibilities;\\
		If $  \iota = \jmath=1, \mu=0 $ then by using Theorm \ref{thm:3.1} $ C_{1,1,0,}$ is an AMDS constacyclic code having the parameters $ [4p^{\varsigma},4p^{\varsigma}-2 ,2].$ Also $ C^{\perp}_{1,1,0} \subseteq C_{1,1,0}.$ Thus, by the above argument, there exists a quantum AMDS code having the parameters $ [[4p^{\varsigma},4p^{\varsigma}-4,2]]_{p^m}.$\\
		If $ \iota=2, \jmath =0, \mu=0 $ then by using Theorm \ref{thm:3.1}   $ C_{2,0,0,}$ is an AMDS constacyclic code having the parameters $ [4p^{\varsigma},4p^{\varsigma}-2 ,2].$ Also $ C^{\perp}_{2,0,0} \subseteq C_{2,0,0}.$ Thus, by the above argument, there exists a quantum AMDS code having the parameters $ [[4p^{\varsigma},4p^{\varsigma}-4,2]]_{p^m}.$
	\end{itemize}
\end{pro}
	Now, under the same procedure as in parts (I) and (III), we obtain the results for the remaining parts.
	\begin{example}
		Let $ p=17, \varsigma=1,m=1.$ Then we have $ x^{68}-1=(x-1)^{17}(x+1)^{17}(x-13)^{17}(x+13)^{17}.$\\	
			Let $ C_{2,0,0,0}= \langle(x-1)^2\rangle.$  It is easy to see that $ C^{\perp}_{2,0,0,0} \subseteq C_{2,0,0,0} .$  By Theorem \ref{thm:4.1}
			is a qAMDS code with  parameters $[[68,64,2]]_{17}$ .\\
	Let $ C_{1,1,0,0}=  \langle(x-1) (x+1)\rangle.$  It is easy to see that $ C^{\perp}_{1,1,0,0} \subseteq C_{1,1,0,0} .$  By Theorem \ref{thm:4.1}  is a code  qAMDS with  parameters $[[68,64,2]]_{17}$.\\
	Let $ C_{2,0,1,0}= \langle(x-1)^2(x-13)\rangle$. It is easy to see that $ C^{\perp}_{2,0,1,0} \subseteq C_{2,0,1,0}.$  By Theorem \ref{thm:4.1}  is a code qAMDS code with  parameters  $ [[68,62,3]]_17  .$
		\end{example}
	
\begin{example}
	Let $ p=7,\varsigma=1,m=1.$ Then we have $ x^{28}-1=(x-1)^7(x+1)^7(x^2+1)^7.$\\ Let $ C_{1,1,0}= \langle(x-1)(x+1)\rangle.$ It is easy to see that $ C^{\perp}_{1,1,0} \subseteq C_{1,1,0} .$  By Theorem \ref{thm:4.1}
	is a  qAMDS code with  parameters $ [[28,24,2]]_{7}$.\\	   

Let $C_{2,0,0}=\langle (x-1)^2\rangle $ It is easy to see that $ C^{\perp}_{2,0,0} \subseteq C_{2,0,0} .$  By Theorem \ref{thm:4.1}
is a qAMDS code with  parameters $[[28,24,2]]_{7}$.	 

\end{example}
\begin{example}
	Let $ p=5,~~~\varsigma=1,~~~m=1.$ Then we have $ x^{20}-2=(x^2+2)^5(x^2+3)^5.$ Let $ C_{0,1}= \langle(x^2-2)\rangle.$ It is easy to see that $ C^{\perp}_{0,1} \subseteq C_{0,1} .$  By Theorem \ref{thm:4.1}
	is a  qAMDS code with  parameters $[[20, 16 ,2]]_{5}.$ 
\end{example}
\begin{example}
Let $ p=3,~~~\varsigma=1,~~~m=1.$ Then we have $ x^{12}-2=(x^2+x+2)^3(x^2+2x+2)^3.$ Let $ C_{1,0}= \langle(x^2+x+2)\rangle.$ It is easy to see that $ C_{1,0}^{\perp}\subseteq C_{1,0}.$ By Theorem \ref{thm:4.1} is a qAMDS code with parameters
  $[[12,8,2]]_{3}$.	
\end{example}

\section{Conclusion}

In this paper, we classify  all AMDS constacyclic codes of length $4p^{\varsigma}$ over the finite field $\mathbb{F}_{p^m}$. To exemplify our findings, we present a range of illustrative examples of AMDS codes specifically constructed over various fields, including $\mathbb{F}_{17}$, $\mathbb{F}_{7}$, $\mathbb{F}_{5}$, and $\mathbb{F}_{3}$.
In addition, we characterize all qAMDS codes that can be constructed from this class of constacyclic codes using the CSS method. This analysis explores the properties and structures of these codes, showcasing how they can be effectively generated within the framework established by constacyclic properties.
Looking ahead, we suggest that future research should focus on investigating qAMDS codes of length $4p^{\varsigma}$ over finite chain rings, as this could yield new insights and applications in coding theory.

\end{document}